\documentclass[draft]{article}
\usepackage{amssymb,amsmath,theorem,a4wide}

\def\sameenum{}

\let\origiff\iff
\ifx\mydefs\undefined\else \fi
\let\mydefs\relax

\ifx\slide\undefined
  \allowdisplaybreaks[2]
\fi

\ifx\loadcyr\undefined\else
\input cyracc.def
\makeatletter
 at 1\@ptsize pt
 at 1\@ptsize pt
 at 1\@ptsize pt
\makeatother

\fi

\def\gobble#1{}
\def\fixsup#1#2{{#1\let\dp\gobble\mathstrut}^#2_}

\def\bme{\hskip.75em\relax}

\def\iff{\quad\text{iff}\quad}

\def\?{\mathbin?}

\let\model\vDash
\let\nmodel\nvDash

\newbox\circlebox
\setbox\circlebox\hbox{$\bigcirc$}
\def\circled#1{%
  \setbox0\hbox to\wd\circlebox{\hss$#1$\hss}\wd0=0pt
  \box0\copy\circlebox}

\let\fii\varphi
\let\tet\vartheta
\let\ep\varepsilon

\ifx\slide\empty
  \def\greek#1{$\expandafter\greeknum\csname c@#1\endcsname$}
\else
  \def\greek#1{$\mathop{\boldsymbol{\expandafter\greeknum\csname c@#1\endcsname}}$}
\fi
\makeatletter
\def\greeknum#1{\ifcase#1\or\alpha\or\beta\or\gamma\or\delta\or\ep
      \or\digamma\or\zeta\or\eta\or\tet\or\iota\else\@ctrerr\fi}
\makeatother

\def\p#1{\langle#1\rangle}

\def\lh#1{\lvert#1\rvert}

\let\bez\smallsetminus

\let\sset\subseteq
\let\nsset\nsubseteq
\let\ssset\subsetneq
\let\Sset\supseteq

\let\onto\twoheadrightarrow

\newcommand\rpair[3][3em]{\mathrel{%
   \begin{matrix}%
     \strut\smash{\xrightonto{\hbox to#1{\hss$#2$\hss}}}\\[-1.7ex]%
     \strut\smash{\xleftembed[\hbox to#1{\hss$#3$\hss}]{}}%
   \end{matrix}}}
\makeatletter
\newcommand\xrightonto[2][]{\ext@arrow 0359\rightontofill{#1}{#2}}
\newcommand\xleftembed[2][]{\ext@arrow 3095\leftembedfill{#1}{#2}}
\def\leftembedfill{\arrowfill@\leftarrow\relbar\hookleftnoarrow}
\def\rightontofill{\arrowfill@\relbar\relbar\onto}
\def\hookleftnoarrow{\DOTSB\relbar\joinrel\rhook}
\makeatother

\def\fl#1{\lfloor#1\rfloor}

\def\cl#1{\lceil#1\rceil}

\mathchardef\#="2023 

\ifx\busspf\undefined\else
\usepackage{bussproofs}
\EnableBpAbbreviations

\fi

\ifx\symlasy\undefined

  \def\centdot#1{{%
    \setbox0\hbox{$\mathop{#1}$}\dimen0 \ht0
    \setbox0\hbox{$#1$}\advance\dimen0 -\ht0
    \setbox2\hbox to\wd0{\hss$\mathop{\cdot}$\hss}\wd2=0pt
    \lower\dimen0\box2\box0 }}
\else
  
  \def\centdot#1{{%
     \setbox0\hbox{$#1$}%
     \raise0.206\ht0\hbox to\wd0{\hss$\cdot$\hss}%
     \kern-\wd0 \box0 }}
\fi

\let\sls|

\def\Up{{\setbox0\hbox{$\uparrow$}%
         \lower\dp0\hbox to\wd0{\hss\vrule width4pt height.4pt\hss}%
         \kern-\wd0\box0}}
\def\UP{{\setbox0\hbox{$\uparrow$}%
         \lower\dp0\hbox to\wd0{\hss\vrule width4pt height.4pt\hss}%
         \kern-\wd0\copy0\kern-\wd0\raise.35ex\box0}}
\def\Down{{\setbox0\hbox{$\downarrow$}%
         \raise\ht0\hbox to\wd0{\hss\vrule width4pt depth.4pt\hss}%
         \kern-\wd0\box0}}

\newif\ifnadm

\def\doadm{\mathrel{%
   \setbox0 \hbox{$\mathop\vdash$}\dimen0 \ht0
   \setbox0 \hbox{$\vdash$}\advance\dimen0 -\ht0
   \vrule width.8\fontdimen8 \textfont3 height\ht0 depth\dp0
   \mkern-1mu
   \lower\dimen0 \hbox{$\vcenter{%
      \ifnadm
        \setbox0 \hbox{$\scriptstyle\sim\mathstrut$}%
        \hbox{\hbox to\wd0{\hss$\scriptstyle/$\hss}\kern-\wd0 \box0 }%
      \else
        \hbox{$\scriptstyle\sim\mathstrut$}%
      \fi}$}}}

\def\nrstyle#1#2#3{%
  \setbox0\hbox{$#1\bigcirc$}%
  \vcenter{\hbox to\wd0{\hss$#2#3$\hss}}%
  \kern-\wd0\box0 }

\DeclareMathOperator\Th{Th}

\ifx\slide\undefined
  
\fi

\def\st{\expandafter\hat}

\def\N{\mathbb N}
\def\Q{\mathbb Q}
\def\Z{\mathbb Z}
\def\RR{\mathbb R}

\mathcode`\*="0203

\mathchardef\mhyphen="2D
\def\cput(#1)#2{\put(#1){\hbox to0pt{\hss#2\hss}}}
\def\txto{${}\to{}$}

\ifx\slide\undefined

\def\noproof{\leavevmode\unskip\bme\vadjust{}\nobreak\hfill$\qed$\par}
\let\qed\Box
\newenvironment{Pf}[1][]
  {\par\noindent\textit{Proof\optpar{#1}:}\bme\ignorespaces}
  {\noproof\pagebreak[2]\vskip\medskipamount\ignorespacesafterend}
\def\optpar#1{\ifx\relax#1\relax\else\ #1\fi}
\def\qedhere{\relax\ifmmode\eqno\qed\expandafter\aftergroup
                   \else\noproof\fi\noqed}
\def\noqed{\let\noproof\relax}

\theoremstyle{plain}
\ifx\shortthm\undefined
\newtheorem{Thm}{Theorem}[section]
\else
\newtheorem{Thm}{Theorem}
\fi
\newtheorem{Prop}[Thm]{Proposition}
\newtheorem{Cor}[Thm]{Corollary}
\newtheorem{Lem}[Thm]{Lemma}

\newtheorem{Obs}[Thm]{Observation}

\ifx\shortthm\undefined
\def\theCl{\arabic{Cl}}
\fi

\theorembodyfont\upshape
\newtheorem{Def}[Thm]{Definition}
\newtheorem{Rem}[Thm]{Remark}
\newtheorem{Exm}[Thm]{Example}

\newenvironment{Pf*}{\let\qed\qedCl\Pf}\endPf

\fi

\ifx\slide\undefined
  \usepackage[reftex]{theoremref}
\fi
\makeatletter
\def\thmref@flush{%
   \ifx\thmref@last\empty\else
      \ifthmref@comma, \thmref@finaltrue\fi \thmref@commatrue
      \thmref@last \ifx\thmref@stack\empty\else s\fi \thmref@num 0
      \let\do\thmref@one \thmref@stack
      \ifcase\thmref@num\or\space and\else\thmref@finaltrue, and\fi
      ~\ref{\thmref@head}\let\thmref@stack\empty\fi}
\def\thmref@one#1{\ifnum\thmref@num>0,\fi
   \space\ref{#1}\advance\thmref@num 1\relax}
\makeatother

\newif\iflinenumbers
\linenumberstrue

{\catcode`\^^I=13 \catcode`\^^M=13
\gdef\doalgo#1#2\end#{\hbox to\hsize{\hss \let^^I\qquad%
  \def\\^^M{\nobreak\hfil\break\vadjust{}\qquad}%
  \fboxsep1em \linenum0 %
  \fbox{\hsize#1\vbox{%
  \everypar{\advance\linenum1 %
      \hbox to1.2em{%
           \hss\iflinenumbers$\scriptstyle\the\linenum$\hskip.6em\fi}}%
  #2}}\hss}\end}}
\newcount\linenum

\def\key{\relax\ifmmode\expandafter\mathbf\else\expandafter\textbf\fi}

\def\allowhyphens{\nobreak\hskip0pt\relax}

\DeclareRobustCommand*\magiclparen{\ifmmode(\else\textup(\allowhyphens\fi}
\DeclareRobustCommand*\magicrparen{\ifmmode)\else\textup)\fi}
\let\lparen=(  \let\rparen=)
\def\magicparon{\catcode`\(\active\catcode`\)\active}
\def\magicparoff{\catcode`\(12 \catcode`\)12 }
\AtBeginDocument{\ifx\ifPreview\iftrue\else\magicparon\fi}
\magicparon
\let (=\magiclparen  \let )=\magicrparen

\ifx\sameenum\undefined
  
  \ifx\enumup\undefined
    
  \else
    
  \fi

\else

\fi

\magicparoff

\mathchardef\comma=\mathcode`\,
{\catcode`\,=\active \gdef,{\comma\penalty\relpenalty}}

\ifx\slide\undefined

\providecommand\dedic{%
  \message{^^JWARNING: embed explicit dedication in the paper!^^J}%
  \thanks{Supported by
    grant IAA100190902 of GA AV \v CR, Center of Excellence CE-ITI under the grant
    P202/12/G061 of GA \v CR, and RVO: 67985840.}}

\author{Emil Je\v r\'abek\dedic\\[\medskipamount]
Institute of Mathematics of the Czech Academy of Sciences\\
\small \v Zitn\'a 25,
115\:67 Praha 1,
Czech Republic,
email: \texttt{jerabek@math.cas.cz}
}

\else\ifx

\author{Emil Je\v r\'abek}
\email{jerabek@math.cas.cz\\[-1em]http://math.cas.cz/\string~jerabek/}
\institution{Institute of Mathematics of the Czech Academy of Sciences, Prague}

\else 

\author[Emil Je\v r\'abek]{Emil Je\v r\'abek\\[\medskipamount]
   \scriptsize\texttt{jerabek@math.cas.cz}\\\texttt{http://math.cas.cz/\string~jerabek/}}
\institute{Institute of Mathematics of the Czech Academy of Sciences, Prague}

\fi\fi

\let\iff\origiff

\DeclareMathOperator\resm{r}

\def\dedic{\thanks{The research leading to these results has received funding from the European Research Council under
the European Union's Seventh Framework Programme (FP7/2007--2013)~/ ERC grant agreement no.~339691. The Institute of
Mathematics of the Czech Academy of Sciences is supported by RVO: 67985840.}}

\title{Rigid models of Presburger arithmetic}

\begin{document}
\maketitle

\begin{abstract}
We present a description of rigid models of Presburger arithmetic (i.e., $\Z$-groups). In particular, we show that
Presburger arithmetic has rigid models of all infinite cardinalities up to the continuum, but no larger.
\end{abstract}

\section{Introduction}\label{sec:introduction}

Starting with classical Galois theory, it has been a common technique in algebra, logic, and other fields of
mathematics to investigate structures by means of their automorphisms. This works well for structures with rich
automorphism groups; in first-order logic, such models abound: for example, all structures have homogeneous (and highly
saturated) elementary extensions, and by the Ehrenfeucht--Mostowski theorem, we can find structures with arbitrarily
complex automorphism groups that in a sense control the structures.

On the opposite side of the spectrum, we have structures with only a few automorphisms, and in particular, \emph{rigid}
structures (i.e., having no automorphisms save the identity). Unlike models with many automorphisms, the existence of
rigid models is a rather mysterious and whimsical property. On the one hand, many theories have no rigid models at all:
this is true of theories as simple as the theory of the infinite set with no further structure, but we can construct
arbitrarily complex such theories just by throwing in a new function symbol denoting a nontrivial automorphism. On the
other hand, many theories of interest do have rigid models, even arbitrarily large.

For example, there are rigid linear
orders (e.g., well orders), and consequently, rigid structures in any class that is ``sufficiently universal'' so that
it can suitably represent all other relational structures (e.g., the theory of graphs). What is more surprising is that
there are rigid dense linear orders or Boolean algebras of any uncountable cardinality (cf.~\cite{she:univ}).

Interestingly, some theories have a limited amount of certain ``trivial'' rigid models, but also many other rigid
models for non-trivial reasons. For instance, any completion of Peano arithmetic has a unique prime model where each
element is definable; these models are clearly rigid. But in fact, Peano arithmetic has many more rigid models: every
model of $\mathrm{PA}$
has a rigid elementary end-extension of the same cardinality (see Kossak and Schmerl~\cite[Thm.~3.3.14]{kos-schm}). For an example
involving a more tame theory: archimedean real-closed fields (or any archimedean ordered fields, for that matter) are rigid for trivial reasons, but as proved by
Shelah~\cite{she:diam-iv}, it also follows from certain combinatorial principles consistent with set theory that there are
large, non-archimedean real-closed fields.

In this note, we will have a look at another tame theory: Presburger arithmetic, or equivalently, the theory of
$\Z$-groups. Models of Presburger arithmetic with rich automorphism groups were studied by Llewellyn-Jones~\cite{llew}. We will
instead have a look at rigid models of the theory; in fact, we will present their complete description. We will see that there
are, in a sense, both ``trivial'' and ``non-trivial'' rigid models, but the amount of non-triviality is limited, which
manifests in a cardinality bound: Presburger arithmetic only has rigid models of sizes up to~$2^{\aleph_0}$. We also discuss
the simpler case of ``unordered Presburger arithmetic'' $\Th(\p{\Z,{+},1})$, which turns out to have only ``trivial''
rigid models.

We note that the tamest theory of ordered abelian groups, namely of the divisible ones, has no rigid models (e.g.,
$x\mapsto2x$ is always an automorphism). The theory of $\Z$-groups, which are just extensions of~$\Z$ by divisible
ordered groups, can thus be thought of as the first nontrivial case in these parts.

\section{Preliminaries}\label{sec:preliminaries}

In this section, we will review the necessary background information about $\Z$-groups, and ordered groups in general.
For references, basic properties of Presburger arithmetic and $\Z$-groups are discussed in Prestel and
Delzell~\cite[\S4.1]{pres-del:math-log}, Marker~\cite[\S3.1]{marker:mod-th}, and Llewellyn-Jones~\cite[\S2--5]{llew};
in a more general set-up, model theory of modules and abelian groups is introduced in Hodges~\cite[\S
A.1--2]{hodges:mod-th}.

In this paper, all groups are commutative, and all ordered structures are linearly ordered.

We assume familiarity with basic properties of ordered groups. We specifically mention that a group can be ordered iff
it is torsion-free; more generally, if $G$ is a torsion-free group, and $H\sset G$ an ordered subgroup, the order
on~$H$ extends to an order on~$G$. A subgroup $H$ of an ordered group~$G$ is \emph{convex} if $-y\le x\le y$ and $y\in
H$ imply $x\in H$. If $H$ is a convex subgroup of~$G$, the order on~$G$ induces an order on $G/H$. Conversely, if $H$
is a subgroup of $G$, $\le_H$ is an order on~$H$, and $\le_{G/H}$ an order on $G/H$, then there exists a unique order
$\le_G$ on $G$ that makes $\p{H,\le_H}$ a convex subgroup of $\p{G,\le_G}$, and that induces the order $\le_{G/H}$ on
$G/H$. Since $\le_G$ is explicitly defined by
\[0\le_Gg\iff 0<_{G/H}g+H\text{ or }(g\in H\text{ and }0\le_Hg),\]
we will call it the \emph{lexicographic order} induced by $\le_H$ and~$\le_{G/H}$, even though strictly speaking we
should reserve this term for the case where $G$ is a direct sum of $H$ and~$G/H$.

Recall that a subgroup $H$ of~$G$ is \emph{pure} if for any integer $n>0$ (w.l.o.g.\ prime), and any element $a\in H$, if $a=nb$
for some $b\in G$, then $a=nb$ for some $b\in H$. If $G$ is torsion free, $b$ is defined uniquely, hence we can restate
the condition as follows: $b\in H$ whenever $nb\in H$ for some integer $n>0$. In particular, if $H$ is a subgroup of a
torsion-free group $G$, then $\{a\in G:\exists n\in\N^{>0}\,na\in H\}$ is the least pure subgroup of $G$ that includes
$H$.

We will use the fact that divisible groups are injective: in particular, a divisible subgroup $H\sset G$ is always a
direct summand of~$G$, and more generally, if $K\sset G$ is a subgroup with trivial intersection with~$H$, we can write
$G=H\oplus K'$ for some $K'\Sset K$.

\emph{Presburger arithmetic} is officially the theory of the monoid $\p{\N,+}$, but this structure is
bi-interpretable with the more convenient ordered group $\p{\Z,+,{\le}}$, and we will work exclusively with the
latter. (Bi-interpretation preserves the automorphism group, among many other things, hence this change is transparent
for our purposes.)

Models of $\Th(\p{\Z,+,{\le}})$ are called \emph{$\Z$-groups}, and they may be characterized as
discretely ordered groups $G$ such that $G/\Z$ is divisible. Here, an ordered group is \emph{discrete} if there exists
a least positive element, conventionally denoted~$1$; the subgroup generated by~$1$, which is the least nontrivial
convex subgroup of~$G$, gives a canonical embedding of the ordered group~$\Z$ in~$G$, and $G/\Z$ denotes the
corresponding quotient group.

The theory of $\Z$-groups admits partial quantifier elimination: any formula $\fii(\vec x)$ is equivalent to a Boolean
combination of linear inequalities $\sum_in_ix_i\ge n$, where $n_i,n\in\Z$, and congruences $x_i\equiv k\pmod m$, where
$0\le k<m\in\N$.

Let $G$ be a $\Z$-group. For any integer $n>0$, the \emph{residue map} $x\mapsto(x\bmod n)$ is a unique group homomorphism
$\resm_n\colon G\to\Z/n\Z$ such that $\resm_n(1)=1$. These maps combine to a unique group homomorphism $\resm\colon
G\to\hat\Z$ such that $\resm(1)=1$, where
\[\hat\Z=\varprojlim_n\Z/n\Z=\prod_{p\text{ prime}}\Z_p\]
is the profinite completion of~$\Z$, and $\Z_p$ denotes the $p$-adic integers.

We will also consider models of the theory $\Th(\p{\Z,+,1})$, which we will call \emph{unordered $\Z$-groups}. If $G$
is a torsion-free group with a distinguished element~$1\ne0$, let $\Z$ denote the subgroup of~$G$ generated by~$1$;
then $G$ is an unordered $\Z$-group iff $G/\Z$ is divisible, and $1$ is not $p$-divisible for any~$p$ (equivalently:
$\Z$ is a pure subgroup of~$G$; equivalently: $G/\Z$ is torsion-free).

If $\p{G,+,{\le}}$ is a $\Z$-group with least positive element~$1$, then $\p{G,+,1}$ is an unordered $\Z$-group, and
the convexity of $\Z\sset G$ ensures that $\le$ induces an order $\le'$ on the quotient group~$G/\Z$. Conversely, if $\p{G,+,1}$
is an unordered $\Z$-group, the torsion-free group $G/\Z$ can be ordered, and for any order on $G/\Z$,
the corresponding lexicographic order on~$G$ that makes $\Z$ convex gives $G$ the structure of a $\Z$-group. That is,
any unordered $\Z$-group expands to a $\Z$-group, and these expansions are in 1--1 correspondence with orders on
$G/\Z$.

Quantifier elimination for unordered $\Z$-groups takes the form that every formula $\fii(\vec x)$ is equivalent to a
Boolean combination of linear equalities $\sum_in_ix_i=n$, and congruences $x_i\equiv k\pmod m$.

The residue maps $\resm_n\colon G\to\Z/n\Z$ and $\resm\colon G\to\hat\Z$ can also be defined for unordered $\Z$-groups
just like before.

Quantifier elimination immediately implies a characterization of elementary substructures: if $\p{G,+,1}$ is an
unordered $\Z$-group, a substructure $H\sset G$ is elementary iff it is itself an unordered $\Z$-group iff it is a pure
subgroup of $G$. We will call such substructures \emph{$\Z$-subgroups of~$G$}. Likewise, if $\p{G,+,{\le}}$ is a
$\Z$-group, its elementary substructures are its $\Z$-subgroups, and these are exactly the substructures containing~$1$
that are $\Z$-groups.

\section{Leibnizian models}\label{sec:leibnizian-models}

What structures can be said to be rigid for trivial reasons? One case that immediately springs to mind are
\emph{pointwise definable models}: i.e., structures $M$ such that every element $a\in M$ is definable in~$M$ without
parameters. Clearly, a complete theory may have only one pointwise definable model up to isomorphism (if at all). In
our case, the unique pointwise definable $\Z$-group is the standard model $\p{\Z,+,{\le}}$. We are, however, more
interested in non-standard (i.e., non-archimedean) examples.

We may, in fact, loosen the condition a bit: in the obvious argument that pointwise definable structures are rigid, we
do not really need that each element be isolated by a single formula from the rest of the model---it is enough if we can
tell apart any pair of distinct elements. Thus, we are led to the class of \emph{pointwise type-definable models}; following
Enayat~\cite{ena:leib}, we will call them more concisely \emph{Leibnizian models}, as they are models that validate
one reading of Leibniz's law of \emph{the identity of indiscernibles}.
\begin{Def}
A structure $M$ is \emph{Leibnizian} if for any $a\ne b\in M$, there exists a formula $\fii(x)$ (without parameters)
such that $M\model\fii(a)$ and $M\nmodel\fii(b)$. In other words, distinct elements have distinct parameter-free
$1$-types.
\end{Def}

For example, a real-closed field is Leibnizian if and only if it is archimedean. (This follows easily from quantifier
elimination: in a non-archimedean real-closed field, any two infinitely large elements have the same type.) Notice that a Leibnizian structure in
a countable language must have cardinality at most~$2^{\aleph_0}$. An elementary substructure of a Leibnizian model is
Leibnizian.
\begin{Obs}\th\label{obs:leib-rigid}
Leibnizian structures are rigid.
\noproof\end{Obs}

Now, how do Leibnizian models of Presburger arithmetic look like? Curiously, the answer does not depend on availability
of the order.

\pagebreak[2]
\begin{Thm}\th\label{thm:leib-zg}
For any $\Z$-group $\p{G,+,{\le}}$, the following are equivalent:
\begin{enumerate}
\item\label{item:1}
$\p{G,+,{\le}}$ is Leibnizian.
\item\label{item:2}
The unordered $\Z$-group $\p{G,+,1}$ is Leibnizian.
\item\label{item:3}
The residue map $\resm\colon G\to\hat\Z$ is an embedding.
\end{enumerate}
\end{Thm}
\begin{Pf}

\ref{item:2}\txto\ref{item:1} is obvious.

\ref{item:3}\txto\ref{item:2}: Let $a\ne b\in G$. By assumption $\resm(a)\ne\resm(b)$, hence $\resm_n(a)\ne\resm_n(b)$
for some $n>0$; thus, the formula $x\equiv k\pmod n$ separates $a$ from~$b$, where $k=\resm_n(a)$.

\ref{item:1}\txto\ref{item:3}: If the group homomorphism $\resm$ is not injective, there is $a\ne0$ such that
$\resm(a)=0$. We claim that $a\ne2a$ have the same types, witnessing that $G$ is not Leibnizian.

This is clear for formulas of the form $x\equiv k\pmod m$, as $a\equiv0\equiv2a\pmod m$ by the choice of~$a$. As for
linear inequalities, they become somewhat degenerate in one variable: $nx\ge m$ is equivalent to $x\ge\cl{m/n}$ if
$n>0$, and to $x\le\fl{m/n}$ if $n<0$. Thus, all nonstandard elements of the same sign (such as $a$ and~$2a$) satisfy
the same inequalities.
\end{Pf}

Thus, any Leibnizian $\Z$-group embeds in~$\hat\Z$; to see that there are many such groups, we only need to observe
that $\hat\Z$ itself already works:

\begin{Lem}\th\label{lem:zhat}
$\hat\Z$ is an unordered $\Z$-group.
\end{Lem}
\begin{Pf}
$\hat\Z$ is torsion-free as each $\Z_p$ is a domain of characteristic~$0$. Notice that for any prime~$p$, an element
$a\in\hat\Z$ is $p$-divisible iff $a\equiv0\pmod p$. On the one hand, this ensures that $1$ is not $p$-divisible; on
the other hand, for any $a\in\hat\Z$, one of $a,a+1,\dots,a+p-1$ is.
\end{Pf}

\begin{Cor}\th\label{cor:leib-iso}
\
\begin{enumerate}
\item Up to isomorphism, Leibnizian unordered $\Z$-groups are exactly the $\Z$-subgroups of~$\hat\Z$.
\item Up to isomorphism, Leibnizian $\Z$-groups are exactly the $\Z$-subgroups of $\p{\hat\Z,{\le}}$, where $\le$ is
a lexicographic order induced by an order on $\hat\Z/\Z$.
\end{enumerate}
In particular, for any $\aleph_0\le\kappa\le2^{\aleph_0}$, there are non-archimedean Leibnizian $\Z$-groups and
unordered $\Z$-groups of cardinality~$\kappa$.
\noproof\end{Cor}

Notice that distinct $\Z$-subgroups of~$\hat\Z$ are non-isomorphic, as they are the images of their own residue maps.

\section{Non-Leibnizian models}\label{sec:non-leibn-models}

As shown by the example of real-closed fields mentioned in Section~\ref{sec:introduction}, non-Leibnizian rigid models
may be more elusive than Leibnizian models. Before we get to them, we first need to know a little about the general
structure of non-Leibnizian models of Presburger arithmetic, or more precisely, of unordered $\Z$-groups.
\begin{Prop}\th\label{prop:dirsum}
Let $G$ be an unordered $\Z$-group. We can write $G=D\oplus L$, where $D=\ker(\resm)$ is a divisible subgroup of~$G$,
and $L$ is a Leibnizian $\Z$-subgroup of~$G$.
\end{Prop}
\begin{Pf}
Put $D=\ker(\resm)$ as indicated. If $a\in D$ and $n\in\N^{>0}$, we have in particular $\resm_n(a)=0$, hence $a=nb$
for some $b\in G$. Then $n\resm(b)=\resm(a)=0$, hence $\resm(b)=0$ as $\hat\Z$ is torsion-free. Thus, $D$ is divisible.

Since $D\sset G$ is divisible, and $\Z\cap D=0$, we have $G=D\oplus L$ for some subgroup $L\sset G$ such that
$\Z\sset L$. Since direct summands are necessarily pure, this makes $L$ a $\Z$-subgroup of~$G$, and it is Leibnizian as
$L\cap\ker(\resm)=0$.
\end{Pf}

\begin{Rem}
Thus, $D$ (or its dimension as a $\Q$-linear space, if we want it numeric) can serve as a measure of
non-Leibnizianity of~$G$.

In a decomposition $G=D\oplus L$ with $D$ divisible and $L$ a Leibnizian $\Z$-subgroup, $D$ is uniquely
determined as $\ker(\resm)$, and also as the largest divisible subgroup of~$G$. On the other hand, $L$ is not unique as
a subgroup of~$G$, though it is of course unique up to isomorphism as $L\simeq G/D\simeq\operatorname{im}(\resm)$.

We stress that the direct sum decomposition in \th\ref{prop:dirsum} only works at the level of groups; if $G$ is
ordered, there is no telling how the orders of $D$ and~$L$ interact.
\end{Rem}

With no order to complicate matters, \th\ref{prop:dirsum} lends itself to an easy description of automorphisms. But
first a bit of notation:
\begin{Def}
If $G$ is an unordered $\Z$-group, and $H\sset G$ a subgroup, we will write $H'$ for $H/(\Z\cap H)$. If moreover
$f\colon H\to K$, $K\sset G$, is a homomorphism such that $f[\Z\cap H]\sset\Z$, let $f'\colon H'\to K'$ be the induced
homomorphism $f'\bigl(a+(\Z\cap H)\bigr)=f(a)+(\Z\cap K)$.
\end{Def}
\begin{Lem}\th\label{lem:aut-unord}
Let $G$ be an unordered $\Z$-group, and $D$ and $L$ as in \th\ref{prop:dirsum}. Since $D\cap\Z=0$, the natural
quotient map is a group isomorphism of $D$ to~$D'$.

Any automorphism $f$ of~$G$ induces an automorphism $f'$ of the group~$G'$ such that $f'(a)-a\in D'$ for all $a\in G'$.
Such an~$f'$ can be represented in a unique way as
\[f'(d+l)=g'(d)+h'(l)+l\qquad(d\in D',l\in L'),\]
where $g$ is an automorphism of the group~$D$, and $h\colon L\to D$ is a group homomorphism with $\Z\sset\ker(h)$.

Conversely, any $f'$ as above is induced by a unique automorphism $f$ of~$G$, namely
\[f(d+l)=g(d)+h(l)+l\qquad(d\in D,l\in L).\]
\end{Lem}

We stress that in the statement above, the language of the structure~$G$ includes the constant~$1$, hence it is fixed
by all automorphism (whereas $f'$ and $g$ are just group automorphisms).

\begin{Pf}
An automorphism $f$ of~$G$ has to fix $\Z$ pointwise, hence it indeed induces a group homomorphism $f'\colon G'\to G'$,
which is in fact an automorphism as $f^{-1}$ induces its inverse.
Moreover, $f$ must preserve congruence formulas, i.e., the residue map $\resm\colon G\to\hat\Z$; thus, $f(a)-a\in
D=\ker(\resm)$ for all $a\in G$, and likewise for~$f'$.

Since $G'=D'\oplus L'$, $f'$ can be uniquely represented as $f'(d+l)=g_1(d)+h_0(l)$, where $g_1\colon D'\to G'$, and
$h_0\colon L'\to G'$. The condition $f'(a)-a\in D'$ further ensures that $g_1\colon D'\to D'$, and that
$h_1(l)=h_0(l)-l$ is a homomorphism $h_1\colon L'\to D'$. If $\pi$ denotes the natural isomorphism of $D$ to~$D'$,
define $g\colon D\to D$ and $h\colon L\to D$ by $g(d)=\pi^{-1}(g_1(\pi(d)))$ and $h=\pi^{-1}(h_1(l+\Z))$. Then $g'=g_1$
and $h'=h_1$. Since $f'$ is an automorphism, so is~$g$.

The converse direction is likewise easy to check.
\end{Pf}
\begin{Rem}\th\label{rem:unord-ord}
If $G$ is a $\Z$-group, then $f\colon G\to G$ is an automorphism of $\p{G,+,{\le}}$ if and only if it is an
automorphism of the unordered $\Z$-group $\p{G,+,1}$, and it is order-preserving. The latter holds if and only if the
induced group automorphism $f'\colon G'\to G'$ is order-preserving, as $\Z\sset G$ is a convex subgroup fixed by~$f$.
\end{Rem}
\begin{Cor}\th\label{cor:unord-rig}
The only rigid unordered $\Z$-groups are the Leibnizian ones.
\end{Cor}
\begin{Pf}
If $G$ is a non-Leibnizian unordered $\Z$-group, the divisible torsion-free group $D$ in the decomposition
$G=D\oplus L$ is nontrivial. But then $D\simeq D'$ has a nontrivial automorphism $g$, for instance $g(x)=2x$. By
\th\ref{lem:aut-unord}, this lifts to a nontrivial automorphism of~$G$, using, e.g., $h=0$.
\end{Pf}

This was a rather anticlimactic answer. However, we will see that the situation for ordered $\Z$-groups is more
interesting: some non-Leibnizian $\Z$-groups are rigid after all. The reason is that some of the automorphisms provided
by \th\ref{lem:aut-unord} are not order-preserving. Even so, rigidity severely constraints $D$ and~$L$, as we
will see in a short while.

We will work for a while with divisible subgroups of~$G'$. Recall that a torsion-free divisible group is nothing else
than a $\Q$-linear space, hence we may treat it with methods from linear algebra. Also notice that a convex subgroup of a
divisible ordered group is divisible.

\begin{Lem}\th\label{lem:d-auto}
Let $G$ be a rigid $\Z$-group, and $D=\ker(\resm)$. Then $G'$ has no proper convex
subgroup~$C$ such that $C+D'=G'$.
\end{Lem}
\begin{Pf}
Assume for contradiction that $C$ is such a subgroup. By linear algebra, $D'$ has a linear subspace~$D_0$ such that
$D'=D_0\oplus(D'\cap C)$, which implies $G'=D_0\oplus C$. Since $C$ is convex, this is not just a direct sum of groups,
but also a lexicographic product of the corresponding orders. Thus,
\[f'(d+c)=2d+c\qquad(d\in D_0,c\in C)\]
defines an order-preserving automorphism of~$G'$, which is nontrivial as $D_0\ne0$. Since $f'(a)-a\in D'$ for all $a\in
G'$, $f'$ lifts to a (still order-preserving) automorphism of~$G$ by \th\ref{lem:aut-unord,rem:unord-ord}.
\end{Pf}
\begin{Lem}\th\label{lem:l-auto}
Let $G$ be a rigid $\Z$-group, and $D=\ker(\resm)$. If $C$ is a convex subgroup of~$G'$
non\-triv\-ially intersecting~$D'$, then $C+D'=G'$.
\end{Lem}
\begin{Pf}
If not, then $G'/(C+D')$ and $C\cap D'$ are nontrivial $\Q$-linear spaces, hence there exists a nonzero group
homomorphism $f_0\colon G'\to C\cap D'$ that vanishes on $C+D'$. Define a group homomorphism $f'\colon G'\to G'$ by
$f'(x)=x+f_0(x)$. Clearly, $f'(x)-x\in D'$; since $f_0[D']=0$, this means $f'$ has an inverse homomorphism $x-f_0(x)$,
hence it is an automorphism of~$G'$. It is also order-preserving: assume that $x\ge0$. If $x>C$, then $f_0(x)\in C$
implies $f'(x)>C$, and a fortiori $f'(x)\ge0$. If $x\in C$, then $f'(x)=x\ge0$. Thus, $f'$ lifts to a nontrivial
automorphism of~$G$ using \th\ref{lem:aut-unord,rem:unord-ord}.
\end{Pf}
\begin{Thm}\th\label{thm:ord-rig-main}
Let $G$ be a $\Z$-group, and write $G=D\oplus L$ with $D$ divisible and $L$ a Leibnizian $\Z$-subgroup of~$G$ as
in \th\ref{prop:dirsum}. If $G$ is rigid, the following conditions hold:
\begin{enumerate}
\item\label{item:6}
$D$ is archimedean.
\item\label{item:5}
$D$ is cofinal in $G$, or trivial.
\item\label{item:4}
$L$ is cofinal in~$G$.
\end{enumerate}
\end{Thm}
\begin{Pf}

\ref{item:6} and~\ref{item:5}: We again identify $D$ with $D'\sset G'$. If $a\in D'$ is positive, let $C$ be the least
convex subgroup of~$G'$ such that $a\in C$. Then $C+D'=G'$ by \th\ref{lem:l-auto}, hence $C=G'$ by \th\ref{lem:d-auto}.
Thus, integer multiples of $a$ are cofinal in~$G'$ (and in~$G$ after lifting back), including~$D'$.

\ref{item:4}: It suffices to show that $L'$ is cofinal in~$G'$. Let $C$ be the least convex subgroup of~$G'$ that
includes~$L'$. Then $C+D'=G'$, hence $C=G'$ by \th\ref{lem:d-auto}.
\end{Pf}

Notice that in \th\ref{thm:ord-rig-main}, $L$ embeds in~$\hat\Z$, and $D$, being archimedean, embeds in~$\RR$:
\begin{Cor}\th\label{cor:card}
All rigid $\Z$-groups have cardinality at most~$2^{\aleph_0}$.
\noproof\end{Cor}

\th\ref{thm:ord-rig-main} is not yet the end of the story, as the given conditions are only necessary, not sufficient.
The missing condition is of a rather different (geometric) flavour, putting restrictions on how $D$ and $L$ sit next to
each other inside~$G$. On the one hand, it does not seem very illuminating, and on the other hand, it tends to be
satisfied for typical examples that one comes up with in practice (as we will see shortly). But anyway, we formulate it here
for completeness.
\begin{Thm}\th\label{thm:ord-rig-mult}
Let $G$ be a $\Z$-group, written as $G=D\oplus L$ as in \th\ref{prop:dirsum}. Assume that $D$ is archimedean, and $D$
and~$L$ are both cofinal in~$G$.

Let $C$ be the largest convex subgroup of~$G$ such that $C\cap D=0$. By assumption,
the ordered group $G/C$ is archimedean, hence there is a homomorphism of ordered groups $\nu\colon G\to\RR$ with
$\ker(\nu)=C$, unique up to multiplication by a positive real constant. Put $G''=\nu[G]$, $D''=\nu[D]$, and
$L''=\nu[L]$.

The following are equivalent:
\begin{enumerate}
\item\label{item:7}
$G$ is rigid.
\item\label{item:8}
The only $\gamma\in\RR^{>0}$ such that $\gamma D''=D''$ and $(\gamma-1)L''\sset D''$ is $\gamma=1$.
\end{enumerate}
\end{Thm}
\begin{Pf}
Notice that since $D\cap C=0$, $\nu$ restricted to~$D$ is an isomorphism of $D$ to~$D''$. Let
$\nu_D^{-1}\colon D''\to D$ denote its inverse.

\ref{item:7}\txto\ref{item:8}: Let $\gamma>0$ be such that $\gamma D''=D''$ and $(\gamma-1)L''\sset D''$, which implies
$(\gamma-1)G''\sset D''$. Define $f_\gamma\colon G\to G$ by
\[f_\gamma(x)=x+\nu_D^{-1}\bigl((\gamma-1)\nu(x)\bigr).\]
Then $f_\gamma$ is a group homomorphism, and $f_\gamma(x)-x\in D$ for each $x\in G$. We also have $\gamma^{-1}D''=D''$
and $(\gamma^{-1}-1)L''=-\gamma^{-1}(\gamma-1)L''\sset-\gamma^{-1}D''=D''$, hence we can apply the same construction to
$\gamma^{-1}$; a straightforward calculation shows that $f_{\gamma^{-1}}$ is the inverse of~$f_\gamma$. Thus,
$f_\gamma$ is a group automorphism.

In order to see that it is order-preserving, let $x\ge0$. Notice that $\nu(f_\gamma(x))=\gamma\nu(x)$. Thus, if $x>C$,
then $f_\gamma(x)>C$, and a fortiori $f_\gamma(x)\ge0$. If $x\in C$, then $\nu(x)=0$, hence $f_\gamma(x)=x\ge0$ as well.

Thus, $f_\gamma$ is an automorphism of~$G$, hence it is the identity. But then $\gamma=1$.

\ref{item:8}\txto\ref{item:7}: Let $f$ be an automorphism of~$G$. By~\th\ref{lem:aut-unord}, $f$ restricts to an
automorphism of~$D\simeq D''\sset\RR$. Pick any nonzero $a\in D$, and put $\gamma=\nu(f(a))/\nu(a)$. Since $f$ is order-preserving, and
$\nu(a)\Q$ is dense in~$\RR$, we must have $\nu(f(x))=\gamma\nu(x)$ for all $x\in D$; in particular, $\gamma D''=D''$. More
generally, the density of $\nu(a)\Q$ (along with convexity of~$C$) implies that for any $x\in G$, $\nu(f(x))=\gamma\nu(x)$, thus
$(\gamma-1)\nu(x)=\nu(f(x)-x)$. Since $f(x)-x\in D$ by \th\ref{lem:aut-unord}, we obtain $(\gamma-1)G''\sset D''$. By
\ref{item:8}, $\gamma=1$. Thus, $f(x)-x\in C\cap D=0$ for all $x\in G$, i.e., $f$ is the identity automorphism.
\end{Pf}
\begin{Rem}\th\label{rem:c-l}
Recall that the subgroup~$L$ in the decomposition $G=D\oplus L$ is not uniquely determined. In the situation
of \th\ref{thm:ord-rig-mult}, we may assume without loss of generality that $L$ is chosen so that it includes the group
$C$ defined in the statement of the theorem, as $C\cap D=0$. This makes the choice somewhat more canonical.
\end{Rem}

Let us mention some convenient sufficient conditions for \th\ref{thm:ord-rig-mult}~\ref{item:8}.
\begin{Lem}\th\label{lem:mult}
Let $D''$ and $L''$ be nonzero $\Q$-subspaces of~$\RR$ such that $D''\cap L''=0$. Then condition~\ref{item:8} from
\th\ref{thm:ord-rig-mult} is satisfied:
\begin{enumerate}
\item\label{item:11}
if $\dim_\Q(D'')$ is finite, or if $\dim_\Q(D'')<\dim_\Q(L'')$;
\item\label{item:10}
if $L''\nsset D''^{-1}D''D''$, or more generally, if there is no $\alpha\in D''^{-1}D''$ such that $L''\oplus D''\sset\alpha D''$;
\item\label{item:9}
if $D''$ is a subfield of~$\RR$, or more generally, if there exists a subfield $K\sset\RR$ such that $D''\sset K\sset
L''\oplus D''$.
\end{enumerate}
\end{Lem}
\begin{Pf}

\ref{item:11}: Either condition ensures $\dim_\Q(D'')<\dim_\Q(L''\oplus D'')$. However, if $\gamma\ne1$ is such that
$(\gamma-1)L''\sset D''=\gamma D''$, then $(\gamma-1)\cdot-$ is a $\Q$-linear embedding of $L''\oplus D''$ in $D''$.

\ref{item:10}: If $\gamma\ne1$ is such that $\gamma D''=D''$, then $\gamma=a/b$ for some nonzero $a,b\in D''$, thus
$\alpha:=(\gamma-1)^{-1}=b/(a-b)\in D''^{-1}D''$, and $L''\oplus D''\sset\alpha D''$.

\ref{item:9}: If $K\ssset L''\oplus D''$, this follows from~\ref{item:10}, as $D''^{-1}D''D''\sset K$. If $K=L''\oplus
D''$, and $\gamma\ne1$ is such that $(\gamma-1)L''\sset D''=\gamma D''$, then 
\[(\gamma-1)K\sset D''\implies\gamma-1\in D''\sset K\implies(\gamma-1)K=K\implies L''\sset K\sset D'',\]
which is a contradiction.
\end{Pf}
\begin{Exm}\th\label{exm:mult}
In order to see that the condition is not automatic, let $L''=(\pi-1)^{-1}\Q$, and
$D''=\sum_{n\in\Z}\pi^n\Q=\Q[\pi,\pi^{-1}]$. Clearly, $L''$ and $D''$ are $\Q$-linear subspaces of~$\RR$, and
$\gamma=\pi$ satisfies $\gamma D''=D''$ and $(\gamma-1)L''\sset D''$. That $L''\cap D''=0$ amounts to
$\pi^n(\pi-1)^{-1}\notin\Q[\pi]$ for any~$n\in\N$, which follows from the transcendence of~$\pi$.
\end{Exm}

While all these classification results are entertaining, we have yet to show the existence of a single non-Leibnizian
rigid $\Z$-group. We are going to remedy this now. First, a simple example.
\begin{Exm}\th\label{exm:exist}
Let $D_0''$ be a proper subfield of~$\RR$ of cardinality~$2^{\aleph_0}$, and $L_0''=\alpha D_0''$ for some
$\alpha\in\RR\bez D_0''$, so that $D_0''\cap L_0''=0$. By \th\ref{lem:mult}, $D_0''$ and $L_0''$ satisfy condition
\ref{item:8} of \th\ref{thm:ord-rig-mult}. Since $\dim_\Q(L_0'')=\dim_\Q(\hat\Z/\Z)$, we may fix a group isomorphism
$\phi\colon\hat\Z/\Z\to L_0''$, and use it to lift the order from $L_0''\sset\RR$ to $\hat\Z/\Z$. This induces a lexicographic order
on~$\hat\Z$ which makes it a Leibnizian $\Z$-group. We put $G=D_0''\oplus\hat\Z$ as a group, and order it
lexicographically so that $G/\Z\simeq D_0''\oplus(\hat\Z/\Z)\simeq D_0''\oplus L_0''\sset\RR$ using~$\phi$. Then $G$ is a non-Leibnizian rigid $\Z$-group by \th\ref{thm:ord-rig-mult}.
\end{Exm}

We can generalize \th\ref{exm:exist} to show that essentially any choice of the various
parameters consistent with \th\ref{thm:ord-rig-main,thm:ord-rig-mult} (with the convention from \th\ref{rem:c-l}) can
be realized.
\begin{Prop}\th\label{prop:exist}
Let
\begin{itemize}
\item $L_0$ be a $\Z$-subgroup of~$\hat\Z$ other than~$\Z$;
\item $C_0$ be a proper $\Z$-subgroup of~$L_0$;
\item $\le_C$ be an order on the group $C_0/\Z$;
\item $D_0''$ and $L_0''$ be trivially intersecting nonzero $\Q$-subspaces of~$\RR$ such that
condition \ref{item:8} from \th\ref{thm:ord-rig-mult} holds (for example, using
\th\ref{lem:mult}), and $\dim_\Q(L_0'')=\dim_\Q(L_0/C_0)$.
\end{itemize}
Then there exists a rigid $\Z$-group $G$ with direct decomposition $G=D\oplus L$ as in \th\ref{prop:dirsum}, and
largest convex subgroup $C$ trivially intersecting~$D$, such that
\begin{itemize}
\item as a group, $L$ is isomorphic to $L_0$;
\item $C$ is isomorphic to $C_0$, equipped with the lexicographic order induced by $\le_C$;
\item the ordered groups $D/C$ and $L/C$ are isomorphic to $D_0''$ and $L_0''$, respectively.
\end{itemize}
\end{Prop}
\begin{Pf}
We endow $C_0$ with the lexicographic order induced by $\le_C$, making $C_0$ a $\Z$-group. Since
$\dim_\Q(L_0'')=\dim_\Q(L_0/C_0)$, we may fix a group isomorphism of $L_0/C_0$ to $L_0''$, which allows us to transfer
the archimedean order from $L_0''$ to $L_0/C_0$. We endow $L_0$ with the lexicographic order induced by the already
constructed orders on $L_0/C_0$ and~$C_0$, making $L_0$ a $\Z$-group with a convex subgroup~$C_0$.

Put $G=L_0\oplus D_0''$ as a group. We equip $G/C_0=(L_0/C_0)\oplus D_0''$ with the order induced by the fixed
isomorphism to $L_0''\oplus D_0''\sset\RR$, and $G$ with the lexicographic order induced by the orders on $G/C_0$
and~$C_0$, making $C_0$ convex. (This is compatible with the previously constructed order on~$L_0$.)

By construction, $G$ is a $\Z$-group satisfying the listed properties, and by \th\ref{thm:ord-rig-mult}, it is rigid.
\end{Pf}
\begin{Cor}\th\label{cor:card-nonleib}
For any infinite cardinals $\kappa,\lambda\le2^{\aleph_0}$, there are non-Leibnizian rigid $\Z$-groups $G$ such that
$\lh D=\kappa$ and $\lh L=\lambda$, where $G=D\oplus L$ are as in \th\ref{prop:dirsum}.
\noproof\end{Cor}
\begin{Rem}\th\label{rem:invar}
If our goal were to describe the composition of non-Leibnizian rigid $\Z$-groups uniquely up to isomorphism, the data
used in \th\ref{prop:exist} would not be quite satisfactory: on the one hand, one such $\Z$-group can be described in
several different ways because $L_0''$ is not uniquely determined; on the other hand, the given data do not describe a
unique $\Z$-group as we left unspecified the isomorphism of $L_0''$ to~$L_0/C_0$.

We may improve the description as
follows: we replace $L_0''$ with $G_0''$ (a $\Q$-subspace of $\RR$ such that $D_0''\ssset G_0''$ and satisfying an
appropriate version of \th\ref{thm:ord-rig-mult}~\ref{item:8}), and a group isomorphism of $L_0/C_0$ to $G_0''/D_0''$.
Then any such data describes a unique non-Leibnizian rigid $\Z$-group up to isomorphism, and conversely, any non-Leibnizian rigid
$\Z$-group is described by almost unique data, the remaining ambiguity being that $G_0''$ (along with its subspace
$D_0''$) is only defined up to multiplication by a real constant.
\end{Rem}

\section*{Acknowledgement}
I wish to thank Ali Enayat, whose question prompted the investigation leading to this paper. I would also like to thank
Roman Kossak for interesting comments, and the anonymous referee for suggestions that helped to improve the clarity of
the paper.

\bibliographystyle{mybib}
\bibliography{mybib}
\end{document}
